\documentclass{amsart}

\usepackage[english]{babel}
\usepackage[all]{xy}
\usepackage{latexsym}
\usepackage{eufrak}
\usepackage{amssymb}
\usepackage{fancyhdr}

\theoremstyle{plain}
\newtheorem{pro}{Proposition}
\newtheorem{thm}[pro]{Theorem}
\newtheorem{lem}[pro]{Lemma}

\theoremstyle{definition}
\newtheorem*{dei}{Definition}

\newtheorem*{Rq}{\sc Remark}
\newtheorem*{Ex}{\sc Example}

\newcommand{\G}{\mathcal{G}}
\newcommand{\N}{\mathcal{N}}

\newcommand{\coker}{\mathop{\mathrm{Coker}}}
\newcommand{\Colim}{\mathop{\mathrm{Colim}}}
\newcommand{\FS}{\mathcal{FS}}
\newcommand{\Sy}{\mathbb{S}}
\newcommand{\CA}{\mathcal{A}}

\newcommand{\Po}{\mathcal{P}}
\newcommand{\F}{\mathcal{F}}

\newcommand{\II}{I}
\newcommand{\Qo}{\mathcal{Q}}

\newenvironment{proo}{\begin{trivlist} \item{\sc {Proof.}}}
  {\hfill $\square$ \end{trivlist}}

\title{\bf Free monoid in monoidal abelian categories}
\author{\bf Bruno Vallette}

\begin{document}

\maketitle

\begin{abstract}
We give an explicit construction of the free monoid in monoidal
abelian categories when the monoidal product does not necessarily
preserve coproducts. Then we apply it to several new monoidal
categories that appeared recently in the theory of Koszul duality
for operads and props. This gives a conceptual explanation of the
form of the free operad, free dioperad and free properad.
\end{abstract}

\tableofcontents

\section*{Introduction}

The construction of the free monoid in monoidal categories is a
general problem that appears in many fields of mathematics. In a
monoidal category with denumerable coproducts, when the monoidal
product preserves coproducts, the free monoid on an object $V$ is
well understood and is given by the words with letters in $V$ (see
\cite{MacLane1} Chapter VII Section $3$ Theorem $2$). In general,
the existence of the free monoid has been established, under some
hypotheses, by M. Barr in \cite{Barr}. When the monoidal product
preserves colimits over the simplicial category, E. Dubuc
described in \cite{Dubuc} a construction for the free monoid. A
general categorical answer was given by G.M. Kelly in \cite{Kelly}
when the monoidal product preserves colimit on one side. Once
again, its construction requires the tensor product to preserve
colimits. The problem is that the monoidal products that appeared
recently in various domains do not share
this general property. \\

In order to study the deformation theory of algebraic structures
like algebras (e.g. associative, commutative, Lie algebras) and
bialgebras (e.g. associative bialgebras, Frobenius bialgebras, Lie
bialgebras, involutive Lie bialgebras), one models the operations
acting on them with operads, properads or props. Like algebras for
operads, it turns out that many types of (bi)algebras can be
defined as particular modules over a monoid in a monoidal
category. Moreover, in simple cases, one does not need the full
machinery of props. For instance, Frobenius bialgebras and Lie
bialgebras can be modelled by \emph{dioperads} (see \cite{Gan})
whereas associative bialgebras and involutive Lie bialgebras
require the notion of \emph{properads} (see \cite{BV}). These two
notions are monoids which generate bigger props but such that the
associated categories of models are the same. The (co)homology
theories and the lax notion ``up to homotopy'' of a particular
type of (bi)algebras are given by the Koszul duality of operads
\cite{GK}, dioperads \cite{Gan} or properads \cite{BV,
MerkulovVallette07}. To generalize Koszul duality theory from
associative algebras \cite{P} to operads, dioperads and properads,
the first step is to extend to notions of bar and cobar
constructions. These constructions are chain complexes whose
underlying space is based on the (co)free
(co)operad (respectively (co)dioperad and (co)properad).\\

This paper was motivated by these new examples of monoidal
structures and the need to make explicit the associated free
monoids for Koszul duality theories. When the monoidal product
preserves coproducts on one side, G.M. Kelly gave a construction
of the free monoid by means of a particular colimit in Equation
$(23.2)$ page $69$ of \cite{Kelly} (see also H.J. Baues, M.
Jibladze, A.Tonks in \cite{BJT} Appendix B and C. Rezk \cite{R}
Appendix A). This construction applies to operads. Since the other
monoidal products considered here do not preserve coproducts
neither general colimits, we need to refine the arguments. For
monoidal abelian categories verifying mild conditions, we produce
a particular colimit which gives a general construction for the
free monoid. Then we apply this result to make explicit the free
monoid in various contexts. The construction of the free properad
is new. For the other examples, the present construction gives a
conceptual explanation for their particular form based on categories of graphs.\\

This paper is organized as follows. In the two first sections, we
fix some conventions and recall the crucial notion of reflexive
coequalizers. The general construction of the free monoid is given
in Section $3$. In Section $4$, we define the notion of
\emph{split analytic functor}, which provides a sufficient
condition to apply the results of the previous section. The last
section is devoted to the description of the free properad, free
$\frac{1}{2}$-prop, free dioperad, free special prop and free
colored operad.

\section{Conventions}

We recall briefly the main definitions used throughout the text.\\

Let $(\CA,\, \boxtimes,\, I)$ be a monoidal abelian category. One
important goal here is to understand the behavior of the monoidal
product with the coproduct of $\CA$. In the sequel, we will denote
the coproduct in an abelian category by
$$\xymatrix@M=5pt@C=14pt{{A\ } \ar@{>->}[r] \ar[rd]_(0.4){f} & A\oplus B \ar[d]^(0.4){f+g}& \ar@{>->}[l] {\ B}\ar[ld]^(0.4){g}  \\
&C. & } $$ When the maps $f$ and $g$ are evident the image of
$f+g$ will be denote by $A+B$.

\begin{dei}[Multiplication functors]
For every object $A$ of $\CA$, we call \emph{left multiplication
functor} by $A$ (respectively \emph{right multiplication functor}), the
functor defined by $L_A \ : \ X \mapsto A\boxtimes X$ (respectively $R_A
\ : \ X \mapsto X\boxtimes A$).
\end{dei}

\begin{dei}[Biadditive monoidal category]
When the left and right multiplication functors are additive for
every object $A$ of $\CA$, the monoidal abelian category $(\CA,\,
\boxtimes,\, I)$ is said to be \emph{biadditive}.
\end{dei}

In a biadditive monoidal category, one knows how to construct
classical objects such as free monoids. When the category is not
biadditive, one can consider the following object to understand
the default of the monoidal product to be additive.

\begin{dei}[Multilinear part]
Let $A$, $B$, $X$ and $Y$ be objects of $\CA$. We call
\emph{multilinear part} in $X$ the cokernel of the morphism
$$ A\boxtimes Y \boxtimes B \xrightarrow{A\boxtimes
i_Y \boxtimes B} A\boxtimes (X\oplus Y)\boxtimes B,$$ which is
denoted by $A\boxtimes (\underline{X} \oplus Y)\boxtimes B$.
\end{dei}

The multilinear part in $X$ is naturally isomorphic to the kernel
of the application
$$A\boxtimes (X\oplus Y)\boxtimes B \xrightarrow{A\boxtimes \pi_Y \boxtimes B}
A\boxtimes Y\boxtimes B,$$ where $\pi_Y$ is the projection
$X\oplus Y \to Y$. The short exact sequence
$$\xymatrix@C=50pt{A\boxtimes (\underline{X} \oplus Y)\boxtimes B\  \ar@{>->}[r] & A\boxtimes (X\oplus Y)\boxtimes B
\ar@{->>}[r]_(0.55){A\boxtimes \pi_Y \boxtimes B} &
\ar@/_1pc/[l]_(0.45){A\boxtimes i_Y \boxtimes B} A\boxtimes
Y\boxtimes B}
$$ splits and we have naturally
$$A\boxtimes (X\oplus Y)\boxtimes B \cong A\boxtimes (\underline{X} \oplus Y)\boxtimes B \oplus
A\boxtimes Y\boxtimes B. $$

A category is biadditive monoidal if and only if one has
$A\boxtimes (\underline{X} \oplus Y)\boxtimes B=A\boxtimes X
\boxtimes B$ for every objects $A$, $B$, $X$ and $Y$.

\begin{dei}[The simplicial categories $\Delta$ and $\Delta_{\rm face}$]
The class of objects of the \emph{simplicial category} $\Delta$ is
the set of finite ordered sets $[n]=\{0<1<\cdots <n\}$, for $n\in
\mathbb{N}$. And the set of morphisms $Hom_{\Delta}([n],\, [m])$
is the set of order-preserving morphisms from $[n]$ to $[m]$.

For $i=0, \ldots,\, n$, one defines the \emph{face map}
$\varepsilon_i\in Hom_{\Delta}([n],\, [n+1])$ by the following
formula
$$\varepsilon_i(j)=\left\{ \begin{array}{ll}
j & \textrm{if}\quad j<i, \\
j+1 & \textrm{if}\quad j\ge i.
\end{array}
\right.$$ The category $\Delta_{\rm face}$ is the subcategory of
$\Delta$ such that the sets of morphisms are reduced to the
compositions of face maps (and the identities $id_{[n]}$).
\end{dei}

\begin{Rq} The category $\Delta_{\textrm{face}}$ is also denoted by
$\Delta^+$ in the literature.
\end{Rq}

\section{Reflexive coequalizers}

In this section, we recall the properties of reflexive
coequalizers which will play a crucial role in the sequel.\\

The multiplication functors $L_A$ and $R_A$ do not necessarily
preserve cokernels, even when they are additive. Nevertheless, for
some monoidal products (for instance the ones treated in the
sequel), the multiplication functors preserve \emph{reflexive
coequalizers}, which is a weaker version of the notion of
cokernel. For more details about reflexive coequalizers, we refer
the reader to the book of P.T. Johnstone \cite{J}.

\begin{dei}[Reflexive coequalizer]
A pair of morphisms $\xymatrix{X_1 \ar@<-0.5ex>[r]_{d_0}
\ar@<0.5ex>[r]^{d_1} & X_0}$ is said to be \emph{reflexive} is
there exists a morphism $s_0\ : \ X_0 \to X_1$ such that $d_0
\circ s_0=d_1 \circ s_0=id_{X_0}$. A coequalizer of a reflexive
pair is called a \emph{reflexive coequalizer}.
\end{dei}

\begin{pro}\label{epi}
\label{reflexive-coequalizers-epi} Let $\Gamma \ : \CA \to \CA$ be
a functor in an abelian category $\CA$. If $\Gamma$ preserves
reflexive coequalizers then it preserves epimorphisms.
\end{pro}

\begin{proo}
Let $\xymatrix{B \ar@{>>}[r]^{\pi} & C}$ be an epimorphism. Since
$\CA$ is an abelian category, $\pi$ is the cokernel of its kernel
$\xymatrix{{A\ } \ar@{^{(}->}[r]^{i}& {B\ } \ar@{>>}[r]^{\pi}&
C}$.

The cokernel $\pi$ can be written as the reflexive coequalizer of
the following pair
$$\xymatrix@C=50pt{A\oplus B\ar@<-1ex>[r]_{d_0}
\ar[r]^{d_1} & \ar@/_1pc/[l]_{s_0}B \ar@{>>}[r]^{\pi}& C,}$$ where
$d_0=i+id_B$, $d_1=id_B$ and $s_0=i_B$.

Since $\Gamma$ preserves reflexive coequalizers, we get that
$\Gamma(\pi)$ is the coequalizer of $(\Gamma(d_0)$, $\,
\Gamma(d_1))$. Therefore, $\Gamma(\pi)$ is an epimorphism.
\end{proo}

A monoidal product is said to \emph{preserve reflexive
coequalizer} if multiplication functors $L_A$ and $R_A$ preserve
reflexive coequalizers for every object $A$ of $\CA$.

\begin{pro}[Lemma 0.17 of \cite{J}]
\label{preservation} Let $(\CA,\, \boxtimes,\  I)$ be a monoidal
abelian category such that the monoidal product preserves
reflexive coequalizers. Let
$$\xymatrix{M_1 \ar@<-1ex>[r]_{d_0} \ar[r]^{d_1}
 &\ar@/_1pc/[l]_{s_0} M_0 \ar@{>>}[r]^{\pi}&
M} \quad \textrm{and} \quad \xymatrix{M_1' \ar@<-1ex>[r]_{d_0'}
\ar[r]^{d_1'}  & \ar@/_1pc/[l]_{s_0'} M_0' \ar@{>>}[r]^{\pi'}& M'}
$$ be two reflexive coequalizers. Then $M \boxtimes M'$ is the
reflexive coequalizer of
$$\xymatrix@C=40pt{M_1\boxtimes M_1' \ar@<-1ex>[r]_{d_0\boxtimes
d_0'} \ar[r]^{d_1\boxtimes d_1'}  &\ar@/_1.5pc/[l]_{s_0\boxtimes
s_0'} M_0\boxtimes M_0' \ar@{>>}[r]^{\pi \boxtimes \pi'}&
M\boxtimes M'}.$$
\end{pro}

In a monoidal category, when the monoidal product preserves
reflexive coequalizers, it shares the following crucial property
with the multilinear part.

\begin{pro}\label{ImSum}
Let $(\CA,\, \boxtimes,\  I)$ be a monoidal abelian category such
that the monoidal product preserves reflexive coequalizers. Let
$V$ and $W$ be two objects of $\CA$ and let $\iota_A : A
\hookrightarrow W$, $\iota_B : B \hookrightarrow W$ be two
sub-objects of $W$. The following two objects are equal in $\CA$
$$\emph{Im}\big(V \boxtimes ( \underline{(A+B)} \oplus W)\big)=
\emph{Im}(V \boxtimes ( \underline{A} \oplus W)) + \emph{Im}(V
\boxtimes ( \underline{B} \oplus W)),$$ where $\emph{Im}$ is
understood to be the image of
$$V \boxtimes ( \underline{A} \oplus W) \hookrightarrow
V \boxtimes ( A \oplus W) \xrightarrow{V\boxtimes(\iota_A + id_W)}
V \boxtimes W.$$
\end{pro}

\begin{proo}
We prove first that the image $\textrm{Im}(V \boxtimes (
\underline{A} \oplus W))$ is equal to the kernel of $V\boxtimes W
\xrightarrow{V\boxtimes \pi} V\boxtimes W/A$, where $\pi$ is the
cokernel of $A \hookrightarrow W$. This cokernel is the following
reflexive coequalizer
$$\xymatrix@C=50pt{A\oplus W\ar@<-1ex>[r]_{d_0}
\ar[r]^{d_1} & \ar@/_1pc/[l]_{s_0}W \ar@{>>}[r]^{\pi}& W/A,}$$
where $d_0=\iota_A+id_W$, $d_1=\pi_W$ and $s_0=i_W$. By the
assumption, it is preserved by left tensoring with $V$
$$\xymatrix@C=65pt{V\boxtimes(A\oplus W) \ar@<-1ex>[r]_(0.55){V\boxtimes(\iota_A+id_W)}
\ar[r]^(0.55){V\boxtimes \pi_W} & \ar@/_2pc/[l]_(0.45){V\boxtimes
i_W} V\boxtimes W \ar@{>>}[r]^{V\boxtimes \pi}& V\boxtimes W/A}.$$
Hence the kernel of $V\boxtimes \pi$ is the image of
$V\boxtimes(\iota_A+id_W) - V\boxtimes \pi_W$. Since $V\boxtimes
(A\oplus W)=V\boxtimes W \oplus V\boxtimes (\underline{A}\oplus
W)$ and since $V\boxtimes(\iota_A+id_W) - V\boxtimes \pi_W$
vanishes on the first component $V \boxtimes W$, the image of
$V\boxtimes(\iota_A+id_W) - V\boxtimes \pi_W$ is given by the
image on the second component $V\boxtimes (\underline{A}\oplus
W)$. The multilinear part $V\boxtimes (\underline{A}\oplus W)$ is
defined as the kernel of $V\boxtimes \pi_W$, so the image of
$V\boxtimes(\iota_A+id_W) - V\boxtimes \pi_W$ on it, is equal to
the image of  $V\boxtimes(\iota_A+id_W)$.

Therefore, the left hand side of the equation is equal to the
kernel of the following reflexive coequalizer
$$\xymatrix@C=75pt{V\boxtimes(A\oplus B\oplus W) \ar@<-1ex>[r]_(0.6){V\boxtimes(\iota_A+\iota_B+id_W)}
\ar[r]^(0.6){V\boxtimes \pi_W} & \ar@/_2pc/[l]_(0.45){V\boxtimes
i_W} V\boxtimes W \ar@{>>}[r]^(.45){V\boxtimes \pi}& V\boxtimes
W/(A+B)},$$ that is the image of
\begin{eqnarray*}
V\boxtimes(\iota_A+\iota_B+id_W) - V\boxtimes \pi_W &=&
V\boxtimes(\iota_A+\iota_B+id_W) - V\boxtimes(\iota_B+id_W) \\&& +
V\boxtimes(\iota_B+id_W) - V\boxtimes \pi_W.
\end{eqnarray*}
If we decompose $V\boxtimes(A\oplus B\oplus W)$ into
$$V\boxtimes(A\oplus B\oplus W)= V\boxtimes(\underline{A}\oplus B\oplus W)\oplus
 V\boxtimes(\underline{B}\oplus W)\oplus V\boxtimes W,$$
we can see that the image of $V\boxtimes(\iota_A+\iota_B+id_W) -
V\boxtimes \pi_W$ on the first component is equal to
$\textrm{Im}(V \boxtimes ( \underline{A} \oplus W))$. On the
second component, it is equal to $\textrm{Im}(V \boxtimes (
\underline{B} \oplus W)) $. And it vanished on the last component,
which concludes the proof.
\end{proo}

As a direct corollary, we get the same formula for a finite number
of sub-objects $A_1, \ldots, A_n$ of $W$, that is
$$\textrm{Im}\big(V \boxtimes ( \underline{(\Sigma_i A_i)} \oplus W)\big)=
\Sigma_i  \ \textrm{Im}(V \boxtimes ( \underline{A_i} \oplus
W)).$$

\section{Construction of the free monoid}
\label{construction}

In this section, we give the construction of the free monoid. We
work in monoidal abelian category $(\CA,\, \boxtimes ,\, I)$ such
that the monoidal product
preserves reflexive coequalizers and sequential colimits. \\

Associated to every object $V$ of $\CA$, we consider the augmented
object $V_+:=I\oplus V$. The injection of $I$ in $V_+$ is denoted
by $\eta \, : \, I \hookrightarrow V_+$ and the projection of
$V_+$ in $I$ is denoted by $\varepsilon \, : \, V_+
\twoheadrightarrow I$. We define $V_n:=(V_+)^{\boxtimes n}$. By
convention, we have $V_0=(V_+)^0=I$. Let $\FS(V)$ denote the
coproduct $\bigoplus_{n\geqslant 0} V_n$.\\

This object is naturally endowed with degeneracy maps
$$\xymatrix@C=55pt{\eta_i \, :\, V_{n}\cong(V_+)^{\boxtimes i}\boxtimes I
\boxtimes (V_+)^{\boxtimes (n-i)} \ar[r]^-{V_{i}\boxtimes \eta
\boxtimes V_{n-i}} & (V_+)^{\boxtimes i}\boxtimes V_+ \boxtimes
(V_+)^{\boxtimes (n-i)}=V_{n+1}}.$$

\begin{Rq}
When $V$ is an augmented monoid, consider the kernel $\bar{V}$ of
the augmentation $V \to I$, called the \emph{augmentation ideal}.
There are face maps on $\FS(\overline{V})$ which define the
categorical (or simplicial) bar construction on $V$ (see
\cite{MacLane1}).
\end{Rq}

In a biadditive monoidal category, the colimit of $\{ V_n \}_n$ on
the small category $\Delta_{\rm face}$ is isomorphic to
$\bigoplus_{n\in \mathbb{N}} V^{\boxtimes n}$, which corresponds
to ``words'' in $V$. In this case, it gives the construction of
the free monoid. In general, the colimit $\Colim_{\Delta_{\rm
face}}V_n$ is not preserved by the monoidal product. Therefore,
one has to consider some quotient of $V_n$
before taking the colimit on $\Delta_{\rm face}$.\\

Denote by $I\boxtimes A \xrightarrow{\lambda_A} A$ and $A\boxtimes
I \xrightarrow{\rho_A} A$ the natural isomorphisms of the monoidal
category $(\CA,\, \boxtimes,\, I)$. We define the morphism $\tau\,
:\, V \to V_2$ by the following composition
$$\xymatrix@C=50pt{ V \ar[r]^(.4){\lambda_V^{-1}\oplus \rho_V^{-1}}& \II\boxtimes
V \oplus V \boxtimes \II \ar[r]^-{\eta \boxtimes i_V - i_V
\boxtimes \eta }& (I\oplus V)\boxtimes (I \oplus V)=V_2},$$ where
$i_V$ is the inclusion $V \hookrightarrow I\oplus V$.

For every $A$ and $B$ two objects of $\CA$, we consider the
``relation'' object
$$ R_{A,B}:=\textrm{Im} \left( A\boxtimes ( \underline{V} \oplus V_2)\boxtimes B
\hookrightarrow
  A\boxtimes (V \oplus V_2)\boxtimes B
\xrightarrow{A\boxtimes (\tau + id_{V_2}) \boxtimes B} A\boxtimes
V_2\boxtimes B \right).$$ We denote by $R_{i,\, n-i-2}$ the
subobject $R_{V_i,\, V_{n-i-2}}$ of $V_n$ for $0\leq i \leq n-2$.

\begin{dei}[$\widetilde{V}_n$]
We define the object $\widetilde{V}_n$ by the formula
$$\widetilde{V}_n:=\coker\left(\bigoplus_{i=0}^{n-2} R_{i,n-i-2}
 \to V_n     \right) .$$
We denote by $R_n=\sum_{i=1}^{n-2}R_{i,n-i-2}$ the image of
$\bigoplus_{i=0}^{n-2} R_{i,n-i-2}$ in $V_n$. Hence, we also
denote $\widetilde{V}_n$ by $V_n/{\textstyle
(\sum_{i=0}^{n-2}R_{i,n-i-2})}=Vn/R_n$ and the following short
exact sequence by
$$\xymatrix{0 \ar[r] & {R_n\ } \ar@{^{(}->}[r]^{i_n} & {V_n\ }
\ar@{>>}[r]^{\pi_n} & \widetilde{V}_n \ar[r] & 0.}$$

\end{dei}

\begin{lem}
$\ $
\begin{enumerate}
\item The morphisms $\eta_i$ between $V_n$ and $V_{n+1}$ induce
morphisms $\widetilde{\eta}_i$ between the quotients
$\widetilde{V}_n$ and $\widetilde{V}_{n+1}$. \item For every
couple $i$, $j$, the morphisms $\widetilde{\eta}_i$ and
$\widetilde{\eta}_j$ are equal.
\end{enumerate}
\end{lem}

\begin{proo}
$\ $
\begin{enumerate}
\item It is enough to see that
$$\left\{ \begin{array}{ll}
\eta_i(R_{j,\, n-j-2})\subset R_{j,\,n-j-1} & \textrm{if} \quad  j\le i-2, \\
\eta_i(R_{j,\, n-j-2})\subset R_{j+1,\,n-j-2} &
\textrm{if}\quad j\ge i,  \\
\eta_i(R_{i-1,\, n-i-1})\subset R_{i,\,n-i-1}+R_{i-1,\,n-i}&
\textrm{for} \quad i=1,\ldots,\, n-1.
\end{array} \right.$$
\item Since $({\eta}_i-{\eta}_{i+1})(V_n)\subset R_{i,\, n-i-1}$,
one has $\widetilde{\eta}_i=\widetilde{\eta}_{i+1}$.
\end{enumerate}
\end{proo}

We denote by $\widetilde{\eta}$ any map $\widetilde{\eta}_i$.

\begin{dei}[$\F(V)$]
The object $\F(V)$ is defined by the following sequential colimit
$$ \xymatrix@M=10pt{\II \ar[r]^-{\widetilde{\eta}} \ar[rd]_(0.35){j_0} &
\widetilde{V}_1=V_1=V_+ \ar[r]^-{\widetilde{\eta}}
 \ar[d]^(0.45){j_1} &
\widetilde{V}_2  \ar[r]^-{\widetilde{\eta}}\ar[dl]_(0.55){j_2}&
\widetilde{V}_3\ar[r]^-{\widetilde{\eta}}
 \ar[dll]_(0.53){j_3}&
\widetilde{V}_4 \,
\cdots \ar[dlll]^(0.35){j_4} \\
& \F(V):=\Colim_{\mathbb{N}} \widetilde{V}_n .& & &}$$
\end{dei}

The colimit $\Colim_{\Delta_{\rm face}} V_n$ has been transformed
to the sequential colimit $\Colim_{\mathbb{N}} \widetilde{V}_n$ by
considering the quotients $\widetilde{V_n}$. The hypothesis that
the monoidal product $\boxtimes$ preserves such colimits gives the
following property.

\begin{lem}
\label{lemColim}
For every object $A$ of $\CA$, the multiplication functors $L_A$ and $R_A$
preserve the previous colimit $\F(V)$. One has
$$A\boxtimes\Colim_{\mathbb{N}} \widetilde{V}_n \cong \Colim_{\mathbb{N}}
 (A\boxtimes \widetilde{V}_n)\quad \textrm{and} \quad
\Colim_{\mathbb{N}} \widetilde{V}_n\boxtimes A \cong
\Colim_{\mathbb{N}} (\widetilde{V}_n \boxtimes A). $$
\end{lem}

We will now endow the object $\F(V)$ with a structure of monoid.\\

The unit $\widetilde{\eta}$ is given by the morphism $ j_0\, :
\, \II \to \F(V)$. The product is defined from the concatenation morphisms
 $V_n\boxtimes V_m \to V_{n+m}$. We consider
$$\mu_{n,\, m} \ : \ \xymatrix{ V_n\boxtimes
V_m \ar[r]^-{\sim}& V_{n+m} \ar@{>>}[r]&
 \widetilde{V}_{n+m} \ar[r]^-{j_{n+m}}& \F(V)}.$$

\begin{pro}
There exists a unique map $\widetilde{\mu}_{n,\, m}\, : \,
\widetilde{V}_n\boxtimes \widetilde{V}_m \to \F(V)$
such that
$$\xymatrix{V_n\boxtimes V_m \ar@{>>}[r]^{\pi_n \boxtimes \pi_m}
\ar[rd]_{\mu_{n,\, m}}&\widetilde{V}_n\boxtimes \widetilde{V}_m
\ar[d]^{\widetilde{\mu}_{n,\, m}} \\
& \F(V).} $$
\end{pro}

\begin{proo}
The cokernels $\widetilde{V}_n$ are reflexive coequalizers of the pairs
$$\xymatrix{R_n\oplus V_n \ar@<-1ex>[r]_(0.55){d_0^n} \ar[r]^(0.55){d_1^n}
  & \ar@/_1pc/[l]_{s_0^n} V_n \ar@{>>}[r]^{\pi_n}&
\widetilde{V}_n,}$$ where $d_0^n=i_n+id_{V_n}$, $d_1^n=\pi_{V_n}$
and
 $s_0^n=i_{V_n}$. Since the monoidal product $\boxtimes$
preserves reflexive coequalizers, we have, by Proposition~\ref{preservation},
that $\pi_n\boxtimes \pi_m$ is the (reflexive) coequalizer of the pair
$(d_0^n \boxtimes d_0^m,\, d_1^n \boxtimes d_1^m)$. Hence, the proof is given
by the universal property of coequalizers. One just has to show that
$\mu_{n,\, m}(d_0^n\boxtimes d_0^m)
=\mu_{n,\, m}(d_1^n \boxtimes d_1^m)$. This relation comes from the
 diagram
$$\xymatrix@C=85pt{(R_n\oplus V_n)\boxtimes(R_m \oplus V_m) \ar[d]^{\pi_{V_n}\boxtimes \pi_{V_m}}
\ar[r]^(0.6){(i_n+id_{V_n})\boxtimes (i_m+id_{V_m})} &
V_n\boxtimes V_m \ar[r]^{\sim} & V_{n+m}
\ar[d]^{\pi_{n+m}}\\
V_n\boxtimes V_m  \ar[r]^{\sim} &V_{n+m} \ar[r]^{\pi_{n+m}} &
\widetilde{V}_{n+m},}$$ which is commutative by the following
arguments. Since $(R_n\oplus V_n)\boxtimes(R_m \oplus V_m)$
decomposes as
\begin{eqnarray*}
(R_n\oplus V_n)\boxtimes(R_m \oplus V_m) &\cong&
(\underline{R_n}\oplus V_n)\boxtimes(R_m \oplus V_m) \oplus
V_n\boxtimes(R_m \oplus V_m)\\
&\cong& \underbrace{(\underline{R_n}\oplus V_n)\boxtimes(R_m
\oplus V_m)}_{(iii)} \oplus
\underbrace{V_n\boxtimes(\underline{R_m} \oplus V_m)}_{(ii)}
\oplus \underbrace{V_n \boxtimes V_m}_{(i)}
\end{eqnarray*}
(see Section $1$), it is enough to prove the relation on each
component.

On $(i)$, we have
\begin{eqnarray*}
&& V_n \boxtimes V_m \hookrightarrow (R_n\oplus V_n)\boxtimes(R_m
\oplus V_m) \xrightarrow{(i_n+id_{V_n})\boxtimes (i_m+id_{V_m})}
V_n \boxtimes
V_m  = \\
&& V_n \boxtimes V_m \hookrightarrow (R_n\oplus V_n)\boxtimes(R_m
\oplus V_m) \xrightarrow{\pi_{V_n}\boxtimes \pi_{V_m}} V_n
\boxtimes V_m = id.
\end{eqnarray*}

On $(ii)$ and on $(iii)$, both composite vanish. The composite
$$V_n\boxtimes(\underline{R_m} \oplus V_m) \hookrightarrow (R_n\oplus V_n)\boxtimes(R_m
\oplus V_m)\xrightarrow{\pi_{V_n}\boxtimes \pi_{V_m}} V_n
\boxtimes V_m$$ is equal to zero by the definition of
$V_n\boxtimes(\underline{R_m} \oplus V_m)$, which is the kernel of
$id_{V_n}\boxtimes \pi_{V_m}$. The other composite
$$V_n\boxtimes(\underline{R_m} \oplus V_m) \hookrightarrow (R_n\oplus V_n)\boxtimes(R_m
\oplus V_m)\xrightarrow{id\boxtimes (i_m+id)   } V_n \boxtimes
V_m\cong V_{n+m} \xrightarrow{\pi_{n+m}} \widetilde{V}_{n+m}$$ is
also equal to zero because the image of
$V_n\boxtimes(\underline{R_m} \oplus V_m) \hookrightarrow
(R_n\oplus V_n)\boxtimes(R_m \oplus V_m)\xrightarrow{id\boxtimes
(i_m+id)   } V_n \boxtimes V_m\cong V_{n+m}$ is a sub-object of
$R_{n+m}=\sum_{j=0}^{m+n-2} R_{j, m+n-j-2}$ by the following
argument. Proposition~\ref{ImSum} shows that this image is equal
to the sum $\sum_{i=0}^{m-2} \textrm{Im}\big( V_n \boxtimes
(\underline{R_{i, m-i-2}} \oplus V_m ) \big)$. Each
$\textrm{Im}\big( V_n \boxtimes (\underline{R_{i, m-i-2}} \oplus
V_m ) \big)$ is a sub-object of the image of the composite
\begin{eqnarray*}
&V_n \boxtimes \big(\underline{V_i \boxtimes (\underline{V}\oplus V_2)
\boxtimes V_{m-i-2}}\oplus V_m) \hookrightarrow  V_n \boxtimes
\big(V_i \boxtimes (V\oplus V_2) \boxtimes V_{m-i-2}\oplus V_m)
\big)&\\& \xrightarrow{V_n\boxtimes( V_i \boxtimes (\tau +
id_{V_2}) \boxtimes V_{m-i-2} +id_{V_m})} V_n \boxtimes V_m,&
\end{eqnarray*}
which is equal to the image of
\begin{eqnarray*}
& V_{n+i} \boxtimes (\underline{V}\oplus V_2) \boxtimes V_{m-i-2}
\hookrightarrow V_{n+i} \boxtimes (V\oplus V_2) \boxtimes
V_{m-i-2}& \\& \xrightarrow{V_{n+i}\boxtimes(  \tau + id_{V_2})
\boxtimes V_{m-i-2} } V_n \boxtimes V_m,&
\end{eqnarray*} that is
$R_{n+i, m-i-2}$.

We apply the same arguments to $(iii)$. Since the image of
$(\underline{R_n}\oplus V_n)\boxtimes(R_m \oplus V_m)
\hookrightarrow (R_n\oplus V_n)\boxtimes(R_m \oplus
V_m)\xrightarrow{(i_n+id)\boxtimes (i_m+id)   } V_n \boxtimes
V_m\cong V_{n+m}$ is a sub-object of the sum
$\sum_{i=0}^{n-2}R_{i,m+n-i-2}$, which is a sub-object of
$R_{n+m}$, the same statement holds for $(iii)$.
\end{proo}

\begin{lem}
There exists a unique morphism $\widetilde{\mu}_{n,\, *}$ such that the
following diagram is commutative
$$\xymatrix@R=30pt@C=40pt{ \widetilde{V}_n \boxtimes \II \ar[r]^-{\widetilde{V}_n\boxtimes
\widetilde{\eta}}
 \ar[d]_(0.57){\widetilde{\mu}_{n,\, 0}}
 \ar[drr] |\hole  &
  \widetilde{V}_n \boxtimes \widetilde{V}_1
  \ar[r]^-{\widetilde{V}_n\boxtimes \widetilde{\eta}}
\ar[dl]_(0.7){\widetilde{\mu}_{n,\, 1}} |(0.25)\hole \ar[dr]
|(.4)\hole & \widetilde{V}_n \boxtimes \widetilde{V}_2 \ \cdots
\ar[dll]^(0.65){\widetilde{\mu}_{n,\, 2}}
\ar[d]\\
\F(V) & &\ar[ll]^-{\exists !\,  \widetilde{\mu}_{n,\, *}}
\widetilde{V}_n \boxtimes \F(V)=\Colim_{\mathbb{N}}(\widetilde{V}_n\boxtimes
\widetilde{V}_*) .}$$
\end{lem}

\begin{proo}
Since the morphisms  $\widetilde{\mu}_{n,\, m}$ commute with the morphisms
$\widetilde{V}_n\boxtimes \widetilde{\eta}$
$$\xymatrix@C=30pt{ \widetilde{V}_n\boxtimes \widetilde{V}_m
\ar[r]^-{\widetilde{V}_n \boxtimes \widetilde{\eta}}
\ar[d]_-{\widetilde{\mu}_{n,\, m}}
& \widetilde{V}_n\boxtimes \widetilde{V}_{m+1} \ar[dl]^-{\widetilde{\mu}_{n,\, m+1}   }\\
 \F(V), & }$$
we have by the universal property of colimits that there exists a unique
map $$\widetilde{\mu}_{n,\, *}\, : \,
\Colim_{\mathbb{N}}(\widetilde{V}_n\boxtimes \widetilde{V}_*) \to
\F(V)$$ such that the diagram commutes. We conclude the proof with
Lemma~\ref{lemColim} which asserts that
$\Colim_{\mathbb{N}}(\widetilde{V}_n\boxtimes
\widetilde{V}_*)=\widetilde{V}_n \boxtimes \F(V)$.
\end{proo}

\begin{lem}
There exists a unique morphism $\bar{\mu}$ such that the following diagram
is commutative
$$\xymatrix@R=30pt@C=40pt{ \II \boxtimes \F(V) \ar[r]^-{\widetilde{\eta}
\boxtimes \F(V)}
 \ar[d]_(0.57){\widetilde{\mu}_{0,\, *}}
 \ar[drr] |\hole  &
 \widetilde{V}_1 \boxtimes \F(V)
   \ar[r]^-{\widetilde{\eta} \boxtimes \F(V)}
\ar[dl]_(0.7){\widetilde{\mu}_{1,\, *}} |(0.25)\hole \ar[dr]
|(.4)\hole & \widetilde{V}_2\boxtimes \F(V) \ \cdots
\ar[dll]^(0.65){\widetilde{\mu}_{2,\, *}}
\ar[d]\\
\F(V) & &\ar[ll]^-{\exists !\,  \bar{\mu}} \F(V) \boxtimes
\F(V)=\Colim_{\mathbb{N}}(\widetilde{V}_n\boxtimes \F(V)) .}$$
\end{lem}

\begin{proo}
The arguments are the same.
\end{proo}

\begin{Rq}
The construction of $\bar{\mu}$ with the first colimit on the left and the
second on the right gives the same morphism.
\end{Rq}

\begin{pro}
The object $\F(V)$ with the multiplication $\bar{\mu}$ and the unit
$\bar{\eta}$ forms a monoid in the monoidal category $(\CA,\, \boxtimes ,\, \II)$.

Moreover, this monoid is augmented. We denote by $\bar{\F}(V)$ its ideal
of augmentation.
\end{pro}

\begin{proo}
The relation satisfied by the unit is obvious. The associativity
of $\bar{\mu}$ comes from the associativity of the maps
$\mu_{n,\,m}$.

The counit map $\varepsilon$ is defined by taking the colimit of the maps
$$ \xymatrix@C=35pt{R_n=\sum_{i=0}^{n-2}R_{i,\, n-2-i} \
\ar[r] & V_n=(V \oplus \II)^{\boxtimes n}\ar@{>>}[r]^-{\pi_n}
\ar[d]_-{\varepsilon^{\boxtimes n}}
 & \widetilde{V}_n \ar@{-->}[dl]^-{\exists !\,
 \widetilde{\varepsilon^{\boxtimes n}}}   \\
& \II^{\boxtimes n}=\II.  &  }$$
\end{proo}

\begin{thm}[Free monoid]
\label{Freemonoid}
In a monoidal abelian category $(\CA,\, \boxtimes,\, I)$ which admits sequential
colimits and such that the monoidal product preserves sequential colimits and
reflexive coequalizers, the monoid  $(\F(V),\,
\bar{\mu},\, \bar{\eta})$ is free on $V$.
\end{thm}

\begin{proo}
The unit of adjunction is defined by
$$u_V\ : \ \xymatrix{V \ \ar@{^{(}->}[r] & V\oplus \II \ar[r]^-{j_1}& \F(V).}$$
For a monoid $(M,\,  \nu,\, \zeta)$, the counit $c_M \, :\, \F(M) \to M$ is given
by the colimit of the following maps $\widetilde{\nu^n}$
$$ \xymatrix@C=35pt{R_n=\sum_{i=0}^{n-2}R_{i,\, n-2-i}\    \ar[r] &
M_n=(M\oplus \II)^{\boxtimes n}\ar@{>>}[r]^-{\pi_n}
\ar[d]_-{\nu^n\circ (M + \zeta)^{\boxtimes n}}
 & \widetilde{M}_n \ar@{-->}[dl]^-{\exists !\, \widetilde{\nu^n}} \\
& M,  &  }$$
where the morphisms $\nu^n \, :\, M^{\boxtimes n} \to M$ represents $n-1$
compositions of the map $\nu$ with itself. The maps $\widetilde{\nu^n}$ are well
defined since $\nu^n\circ (M
+ \zeta)^{\Box n} (R_{i,\, n-2-i})=0$, for every  $i$, by associativity of $\nu$.

One has immediately the two relations of adjunction
\begin{eqnarray*}
\xymatrix{ \F(V) \ar[r]^-{\F(u_V)}& \F(\F(V))\ar[r]^-{c_{\F(V)}}
&\F(V)}&= &id_{\F(V)}
\quad \textrm{and} \\
\xymatrix{ M \ar[r]^-{u_M}& \F(M) \ar[r]^-{c_{\F(M)}}& M }& =&
id_M.
\end{eqnarray*}
\end{proo}

\section{Split analytic functors}

In this section, we define the notion of \emph{split analytic functor}.
We show that a split analytic functor preserves reflexive coequalizers.
We will use this proposition in the next section to show that the monoidal
products, considered in the sequel, preserve reflexive coequalizers.\\

Let $\CA$ be an abelian category. And denote by
 $\Delta_n$ the diagonal functor $\CA \to \CA^{\times n}$.

\begin{dei}[Homogenous polynomial functors]
We call a \emph{homogenous polynomial functor of degree $n$} any
functor $f \, :\, \CA \to \CA$ that can be written $f=f_n\circ
\Delta_n$ with $f_n$ a  functor $\CA^{\times n} \to \CA$ additive
in each input.
\end{dei}

\begin{dei}[Split polynomial functor]
A functor $f\,  : \, \CA
\to \CA$ is called \emph{split polynomial} if it is the direct sum
of homogenous polynomial functors
$f=\bigoplus_{n=0}^N f_{(n)}$.
\end{dei}

The functors induced by monoidal products can not always be written with
a finite sum of polynomial functors.

\begin{dei}[Split analytic functors]
We call \emph{split analytic functor} any functor
$f\, : \, \CA \to \CA$ equal to $f=\bigoplus_{n=0}^\infty f_{(n)}$ where
$f_{(n)}$ is an homogenous polynomial functor of degree $n$.
\end{dei}

\begin{Ex} The Schur functor $\mathcal{S}_\Po \, : \, \textrm{Vect} \to \textrm{Vect} $
associated to an $\Sy$-module $\Po$ (a collection
$\{\Po(n)\}_{n\in \mathbb{N}}$
of modules over the symmetric groups $\Sy_n$)
defined by the following formula
$$\mathcal{S}_\Po (V) :=\bigoplus_{n=0}^\infty \Po(n)\otimes_{\Sy_n}
V^{\otimes n}$$
is a split analytic functor.
\end{Ex}

\begin{pro}
\label{analyticreflexive}
Let $f=\bigoplus_{n=0}^N f_{(n)}$ be a split analytic functor such that for every $n\in \mathbb{N}$, every $i\in [n]$ and every $X_1, \ldots, X_{i-1}, X_{i+1}, \ldots, X_n \in \CA$ the functor
$X\mapsto f_n(X_1, \ldots, X_{i-1}, X, X_{i+1}, \ldots, X_n)$ preserves reflexive coequalizers. Then $f$ preserves reflexive equalizers.
\end{pro}

\begin{proo}
Let  $\xymatrix{X_1 \ar@<-1ex>[r]_{d_0} \ar[r]^{d_1} &
\ar@/_1pc/[l]_{s_0}X_0 \ar@{>>}[r]^{\pi}& X}$ be a reflexive coequalizer. The result comes from the
the formula
$$\sum_{i=1}^n f_n(X_0,\ldots,\,
\underbrace{(d_0-d_1)(X_1)}_{i^{\textrm{th}}\ \textrm{place}},
\,\ldots ,\ X_0)=\left( f_n(d_0,\ldots,\ d_0)-f_n(d_1,\ldots,\,
d_1)\right)\circ \Delta_n(X_1).$$ The inclusion $\supset$ is always true since
$$f_n(d_0,\ldots,\
d_0)-f_n(d_1,\ldots,\, d_1)=\sum_{i=1}^n f_n(d_0,\ldots ,\, d_0
,\, \underbrace{d_0-d_1}_{i^{\textrm{th}}\ \textrm{place}},\,
d_1,\ldots ,\,d_1).$$

The reverse inclusion $\subset$ lies on $s_0$ and comes from
\begin{eqnarray*}
&& f_n(X_0,\ldots,\,X_0,\,  (d_0-d_1)(X_1),\, X_0 ,\,\ldots ,\
X_0)\\
&=&  f_n(X_0,\ldots,\,  d_0(X_1) ,\,\ldots ,\
X_0)-f_n(X_0,\ldots,\,
d_1(X_1) ,\,\ldots ,\ X_0)\\
&=& f_n(d_0s_0(X_0),\ldots,\,  d_0(X_1) ,\,\ldots ,\
d_0s_0(X_0))-\\
&&f_n(d_1s_0(X_0),\ldots,\, d_1(X_1) ,\,\ldots ,\
d_1s_0(X_0)).
\end{eqnarray*}
\end{proo}

\section{Applications}

The aim of this section is to apply the previous construction of
the free monoid of new families of monoidal categories that
appeared recently in the theory of Koszul duality. In order to
understand the deformations of algebraic structures, one models
them with an algebraic object (e.g. operads, colored operads,
properads). This algebraic object turns out to be a monoid in an
appropriate monoidal category. The best example is the notion of
operad which is a monoid in the monoidal category of $\Sy$-modules
with the composition product $\circ$ (see J.-L. Loday \cite{Loday}
or J.P. May \cite{May})).

The example of the free properad is new. The other free monoids
given here were already known but the construction of
Section~\ref{construction} gives a conceptual explanation for
their particular form.

\subsection{Free properad}

We recall the definition of the monoidal category of
$\Sy$-bimodules with the connected composition product
$\boxtimes_c$. For a full treatment of $\Sy$-bimodules and the
related monoidal categories, we refer the reader to \cite{BV}.

\begin{dei}[$\Sy$-bimodules]
An \emph{$\Sy$-bimodule} $\Po$ is a collection $(\Po(m,\,
n))_{m,\, n\in \mathbb{N}}$ of modules over the symmetric groups
$\Sy_m$ on the left and $\Sy_n$ on the right, such that the two
actions are compatible. We denote the category of $\Sy$-bimodules
by $\Sy$-biMod.
\end{dei}

An $\Sy$-bimodule models the operations with $n$ inputs and $m$
outputs acting on a type of algebraic structures (like algebras,
bialgebras for instance).

In order to represent the possible compositions of these
operations, we introduced in \cite{BV} a monoidal product
$\boxtimes_c$ in the category of $\Sy$-bimodules. The product $\Qo
\boxtimes_c \Po$ of two $\Sy$-bimodules is given by the sum on
connected directed graphs with $2$ levels $\mathcal{G}^2_c$ where
the vertices of the first level are indexed by elements of $\Po$
and the vertices of the second level are indexed by elements of
$\Qo$ (see Figure~\ref{produit}).

\begin{figure}[h]
$$ \xymatrix{
\ar[dr]_(0.7){1} 1&2 \ar[d]_(0.6){2} &4\ar[dl]^(0.7){3} &3\ar[dr]_(0.7){1} & &5\ar[dl]^(0.7){2} \\
\ar@{--}[r]& *+[F-,]{p_1} \ar@{-}[dl]_(0.3){1}
\ar@{-}[drr]^(0.3){2} \ar@{--}[rrr]& & & *+[F-,]{p_2}
\ar@{-}[dr]^(0.3){2}   \ar@{-}[dll]_(0.3){1} |(0.75) \hole  \ar@{--}[r]& \\
 *=0{} \ar[dr]_(0.7){1} & &*=0{}\ar[dl]^(0.7){2} &*=0{}\ar[dr]_(0.7){1} & &*=0{} \ar[dl]^(0.7){2}\\
\ar@{--}[r] & *+[F-,]{q_1} \ar[d]_(0.3){1}  \ar@{--}[rrr]& & &
*+[F-,]{q_2} \ar[dl]_(0.3){1} \ar[d]_(0.3){2} \ar[dr]^(0.3){3}
\ar@{--}[r] & \\
&4 & &1 &2 &3 } $$ \caption{Example of an element of $\Qo \boxtimes_c \Po$.}
\label{produit}
\end{figure}

We denote by $In(\nu)$ and
$Out(\nu)$ the sets of inputs and outputs of a vertex $\nu$ of a
graph. Let $\mathcal{N}_i$ be the set of vertices on the $i^{\textrm{th}}$ level.

\begin{dei}[Connected composition product $\boxtimes_c$]
Given two $\Sy$-bimodules $\Po$ and $\Qo$, we define their product
by the following formula

$$\Qo\boxtimes_c \Po := \left( \bigoplus_{g\in\G^2_c}
\bigotimes_{\nu \in \N_2} \Qo(|Out(\nu)|,\, |In(\nu)|) \otimes_k
\bigotimes_{\nu \in \N_1} \Po(|Out(\nu)|,\, |In(\nu)|)\right)
\Bigg/ \approx \ ,$$ where the equivalence relation $\approx$ is generated
by
$$
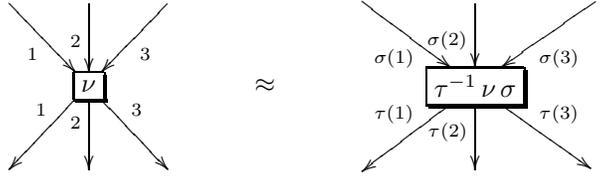
\xymatrix{\ar[dr]_(0.5){1} &\ar[d]_(0.5){2} & \ar[dl]^(0.5){3}
& & \ar[dr]_(0.5){\sigma(1)}& \ar[d]_(0.5){\sigma(2)}&
\ar[dl]^(0.5)
{\sigma(3)} \\
 & *+[F-,]{\nu} \ar[dl]_(0.4){1}\ar[d]_(0.4){2}\ar[dr]^(0.4){3}&
 &\approx & &
 *+[F-,]{\tau^{-1}\, \nu\,\sigma} \ar[dl]_(0.5){\tau(1)}\ar[d]_(0.5){\tau(2)}
 \ar[dr]^(0.5){\tau(3)}&\\
 & & & & & &\ .}$$
\end{dei}

The $\Sy$-bimodule $I$ defined by the following formula
$$I:=\left\{
\begin{array}{l}
I(1,\, 1)=k,    \\
I(m,\,n)=0 \quad \textrm{otherwise}.
\end{array} \right.$$
plays the role of the unit in the monoidal category $(\Sy\textrm{-bimod},\,
\boxtimes_c)$. It corresponds to the identity operation.

\begin{dei}[Properads]
We call a \emph{properad} a monoid $(\Po,\, \mu,\, \eta)$ in the monoidal
category $(\Sy\textrm{-bimod},\,  \boxtimes_c,\, I)$.
\end{dei}

We are going to describe the free properad. To do that, we first show the following
lemma.

\begin{lem}
For every pair $(A,\, B)$ of $\Sy$-bimodules, the functor
$$\Phi_{A,\, B}\, :\, X \mapsto
A\boxtimes_c X \boxtimes_c B$$ is a split analytic functor.
\end{lem}

\begin{proo}
The $\Sy$-bimodule $A\boxtimes_c X \boxtimes_c B$ is given by the
direct sum on $3$-level connected graphs $\mathcal{G}^3_c$ such that the
vertices of the first level are indexed by elements of $B$, the
vertices of the second level are indexed by elements of $X$ and
the vertices of the third level are indexed by elements of $A$.
Denote by $\mathcal{G}^3_{c,\, n}$ the set of 3-level graphs with $n$
vertices on the second level. Therefore, the functor $\Phi_{A,\,
B}$ can be written
\begin{eqnarray*}
\Phi_{A,\, B}(X) &=& A\boxtimes_c X \boxtimes_c B \\
&=& \bigoplus_{n\in \mathbb{N}} \Big( \bigoplus_{g \in
\mathcal{G}_{c,\, n}^3} \bigotimes_{\nu \in \mathcal{N}_1}
A(|Out(\nu)|,\, |In(\nu)|) \otimes \bigotimes_{i=1}^n
X(|Out(\nu_i)|,\, |In(\nu_i)|) \otimes\\
&& \bigotimes_{\nu \in \mathcal{N}_3}
B(|Out(\nu)|,\, |In(\nu)|) \Big) \Big/ \approx  \\
&=& \bigoplus_{n\in \mathbb{N}} {\Phi_n}(X, \ldots, X),
\end{eqnarray*}
where $\Phi_n$ is an homogenous polynomial functor of degree $n$.
\end{proo}

\begin{pro}
\label{SbiModMultilineaire} The category
$(\Sy\textrm{-biMod},\, \boxtimes_c,\, I)$ is a
monoidal abelian category that preserves reflexive coequalizers
and sequential colimits.
\end{pro}

\begin{proo}
For every $\Sy$-bimodule $A$, the left and right multiplicative
functors $L_A:=A\boxtimes_c \bullet$ and $R_A:=\bullet \boxtimes_c
A$ by $A$ are split analytic functors by the previous lemma.
Since the functors $\Phi_n$ preserve reflexive equalizers in each variable, they preserve reflexive coequalizers by Proposition~\ref{analyticreflexive}.
\end{proo}

This proposition allows us to apply Theorem~\ref{Freemonoid}.
Let us interpret this construction in the framework of $\Sy$-bimodules.

\begin{thm}
\label{proplibre}
The free properad on an $\Sy$-bimodule $V$ is given by the sum
on connected graphs (without level) $\G$ with the vertices indexed by elements
of $V$
$$\F(V)=\left( \bigoplus_{g\in \mathcal{G}_c} \bigotimes_{\nu \in
\mathcal{N}} V(|Out(\nu)|,\, |In(\nu)|)\right) \Bigg/ \approx \
.$$
The composition $\mu$ comes from the composition of directed graphs.
\end{thm}

\begin{proo}
The multilinear part in $Y$, denoted $A\boxtimes_c (X\oplus
\underline{Y}) \boxtimes_c B$ is isomorphic to the sub-$\Sy$-bimodule
of $A\boxtimes_c (X\oplus Y) \boxtimes_c B$ composed by 3-level connected graphs
with the vertices of the second level indexed by elements of $X$ and
at least one element of $Y$. Let $V$ be an
$\Sy$-bimodule. Denote by $V_+=I\oplus V$
the augmented $\Sy$-bimodule. Consider the $\Sy$-bimodule
$V_n:=(V_+)^{\boxtimes_c n}$ given by $n$-level connected graphs where the vertices are
indexed by elements of $V$ and $I$. The $\Sy$-bimodule
$\widetilde{V_n}:=\coker \left( \bigoplus_i
R_{V_i,\, V_{n-i-2}}\to V_n \right)$ corresponds to the quotient
of the $\Sy$-bimodule of $n$-level connected graphs by the relation $V\boxtimes_c
I  \simeq I\boxtimes_c V$, which is equivalent to forget the
levels.
\end{proo}

The notion of properad is a ``connected'' version of the notion of
prop (see F.W. Lawvere \cite{Lawvere}, S. Mac Lane \cite{McLane2}
and J.F. Adams \cite{A}). For more details about the link between
these two notions we refer the reader to \cite{BV}. From the
previous theorem, one can get the description of the free prop on
an $\Sy$-bimodule $V$. We find the same construction of the free
prop as B. Enriquez and P.
Etingof in \cite{EE} in terms of forests of graphs without levels.\\

Recall that we have the following inclusions of monoidal abelian
categories (see \cite{BV} Section $1$)
$$(Vect,\, \otimes_k,\, k) \hookrightarrow (\Sy\textrm{-Mod}, \circ,\, I)
\hookrightarrow (\Sy\textrm{-biMod}, \boxtimes_c ,\, I),$$ where
the product $\circ$ of $\Sy$-modules corresponds to the
composition of the related Schur functors. It can be represented
by trees with $2$ levels (see J.-L. Loday \cite{Loday} and J.P.
May \cite{May}). A monoid for the product $\circ$ is called an
\emph{operad}. A direct corollary of the preceding theorem gives
the free associative algebra and the free operad as the direct sum
on trees without levels. Since the monoidal product $\circ$ of
$\Sy$-modules preserves coproducts on the left, the free operad
can be given by more simple colimit (see Kelly \cite{Kelly}
Equation $(23.2)$ page $69$,  Baues-Jibladze-Tonks \cite{BJT}
Appendix B and Rezk \cite{R} Appendix A).

\subsection{Free $\frac{1}{2}$-prop}

On the category of $\Sy$-bimodules, one can define three other
monoidal products. When one wants to model the operations acting
on types of (bi)algebras defined by relations written with simple
graphs (without loops for instance), there is no need to use the
whole machinery of properads. It is enough to restrict to simpler
types of compositions. That is we  consider monoidal products
based on these compositions. The main property is that the
category of (bi)algebras over this more simple object is equal to
the category of (bi)algebras of the associated properad.
Therefore, in order to study the deformation theory of these
(bi)algebras, it
is enough to prove Koszul duality theory for the simpler monoid. (For more details on these
notions, we refer the reader to the survey of M. Markl \cite{Markl2}).\\

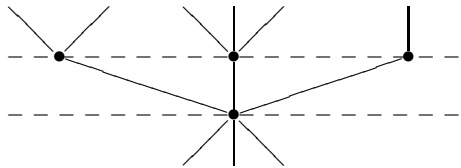
\begin{figure}[h]
$$ \xymatrix@R=16pt@C=16pt{ \ar@{-}[dr]& & \ar@{-}[dl]&  \ar@{-}[dr] &\ar@{-}[d] &\ar@{-}[dl] & & \ar@{-}[d]   & \\
\ar@{--}[rrrrrrrr]& *{\bullet}\ar@{-}[drrr]& & & *{\bullet} \ar@{-}[d] & & &*{\bullet} \ar@{-}[dlll] & \\
\ar@{--}[rrrrrrrr]& & & & *{\bullet}  \ar@{-}[dr]\ar@{-}[d]\ar@{-}[dl]& & & & \\
& & & & & & & & } $$ \caption{Example of a graph in
$\mathcal{G}^{\frac{1}{2}}_2$.} \label{halfGraph}
\end{figure}

Denote by $\mathcal{G}^{\frac{1}{2}}_2$ the set of $2$-level
connected graphs such that every vertices of the first level has
only one output or such that every vertices of the second level
has only one input (see Figure~\ref{halfGraph}).

\begin{dei}[Product $\Box_{\frac{1}{2}}$]
Let $\Po$, $\Qo$ be two $\Sy$-bimodules. Their product $\Qo
\Box_{\frac{1}{2}} \Po$ is the restriction of the connected
composition product $\Qo \boxtimes_c\Po$ on graphs of
$\mathcal{G}^{\frac{1}{2}}_2$.
\end{dei}

\begin{figure}[h]
$$ \xymatrix@R=12pt@C=12pt{ \ar@{-}[dr]& & \ar@{-}[dl]& & & & & \\
 & *{\bullet}\ar@{-}[dr] & &\ar@{-}[dl] & \ar@{-}[dd]&\ar@{-}[dr] & &\ar@{-}[dl] \\
 & & *{\bullet} \ar@{-}[drr]& & & &*{\bullet} \ar@{-}[dll]& \\
 & & & &*{\bullet} \ar@{-}[dr]\ar@{-}[dll]& & & \\
 & & *{\bullet}\ar@{-}[dl]\ar@{-}[d]\ar@{-}[dr]& & &*{\bullet}\ar@{-}[dl]\ar@{-}[dr] & & \\
 & *{\bullet}\ar@{-}[dl]\ar@{-}[dr] & & & & & & \\
 & & & & & & & } $$ \caption{Example of the underlying graph in a free $\frac{1}{2}$-
PROP.} \label{freehalf}
\end{figure}
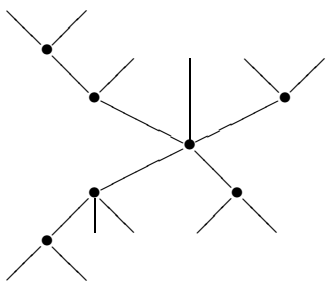

This product is associative and has $I$ for unit. Therefore,
$(\Sy\textrm{-biMod},\, \Box_{\frac{1}{2}}, \, I)$ is a monoidal
abelian category. A monoid in this category is a
$\frac{1}{2}$-prop, notion defined by M. Markl and A.A. Voronov in
\cite{MV} and introduced by M. Kontsevich \cite{K, Markl}. Once
again, we can apply Theorem~\ref{Freemonoid}. The free
$\frac{1}{2}$-prop on an $\Sy$-bimodule $V$ is given the sum on
graphs with one vertex in the middle, grafted above by trees
(without levels) and grafted below by reversed trees without
levels (see Figure~\ref{freehalf}).

\subsection{Free dioperad}

\begin{figure}[h]
$$ \xymatrix@R=16pt@C=16pt{\ar@{-}[dr] & \ar@{-}[d]& \ar@{-}[dl]&\ar@{-}[dr] &  &\ar@{-}[dl] \\
\ar@{--}[rrrrr]&*{\bullet}\ar@{-}[d]\ar@{-}[drrr]  & & &*{\bullet}\ar@{-}[d] & \\
\ar@{--}[rrrrr]& *{\bullet} \ar@{-}[dr]\ar@{-}[d]\ar@{-}[dl]& & & *{\bullet}\ar@{-}[dl]\ar@{-}[dr]& \\
& & & & & } $$ \caption{Example of a connected graph with $2$
levels.} \label{genus0}
\end{figure}
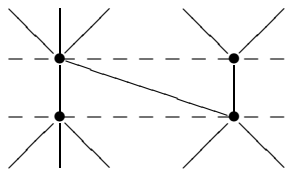

Wee Liang Gan in \cite{Gan} considered the case when the permitted
compositions are based on graphs of genus $0$ (see
Figure~\ref{genus0}).

\begin{dei}[Product $\Box$]
The product $\Qo \Box \Po$ of two $\Sy$-bimodules is given by the restriction
on $2$-level connected graphs of genus $0$ of the monoidal product $\boxtimes_c$.
\end{dei}

Once again, this defines a new monoidal category structure on
$\Sy$-bimodules. A monoid for this product corresponds to the
notion of \emph{dioperad} introduced in \cite{Gan}. By the same
arguments, Theorem~\ref{Freemonoid} shows that the free dioperad
on an $\Sy$-bimodule $V$ is given by the direct sum of graphs of
genus $0$, without levels, whose vertices are
indexed by elements of $V$.  \\

For example, Lie bialgebras, Frobenuis algebras, infinitesimal
bialgebras can be modelled by a dioperad
(see \cite{Gan}). \\

\subsection{Free special prop}

In order to give the resolution of the prop of bialgebras, M.
Markl in \cite{Markl} defined the notion of \emph{special props}.
It is corresponds to monoids in the monoidal category of
$\Sy$-bimodules where the monoidal product is based only on
composition called \emph{fractions} (\cite{Markl} definition
$19$). We can apply Theorem~\ref{Freemonoid} in this case which
gives the free special prop.

Notice that this notion of special props corresponds to the notion
of \emph{matrons} defined by S. Saneblidze R. Umble in \cite{SU}
and is related to the notion of $\frac{2}{3}$-prop of B. Shoikhet
\cite{Sh}.

\subsection{Free colored operad}

Roughly speaking, a colored operad is an operad where the
operations have colors indexing the leaves and the root. The
composition of such operations is null if the colors of the roots
of the inputs operations do not fit with the colors of the
operation below. C. Berger and I. Moerdijk defined a monoidal
product of the category of colored collections such that the
related monoids are exactly colored operad (see Appendix of
\cite{BM}). Once again, Theorem~\ref{Freemonoid} applies in this
case and we get the description of the free colored operad by
means of trees without levels.

\bigskip

\begin{center}
\textsc{Acknowledgements}
\end{center}

I would like to thank Michael Batanin, Benoit Fresse, Jean-Louis
Loday and Ross Street for useful discussions, advice and
improvements of this paper. I am very grateful to the referee for his
numerous remarks and for his knowledge of Kelly's papers.

\bigskip


\bigskip

{\small \textsc{Laboratoire J.A. Dieudonn\'e, Universit\'e de Nice
Sophia-Antipolis, Parc Valrose, 06108 Nice
Cedex 02, France}\\
E-mail address : \texttt{brunov@math.unice.fr}\\
URL : \texttt{http://math.unice.fr/$\sim$brunov}}

\end{document}